\documentclass{amsart}

\usepackage{amssymb}

\newtheorem{thm}{Theorem}

\newtheorem{lem}[thm]{Lemma}

\newtheorem{prop}[thm]{Proposition}

\theoremstyle{definition}
\newtheorem{defn}[thm]{Definition}

\newtheorem{say}[thm]{}

\newtheorem{ques}[thm]{Question}    

\newtheorem{rem}[thm]{Remark}          

\newtheorem*{ack}{Acknowledgments}      

\newtheorem{defn-thm}[thm]{Definition--Theorem}  
\newtheorem{defn-lem}[thm]{Definition--Lemma}  

\theoremstyle{remark}

\renewcommand{\c}[0]{{\mathbb C}}  

\renewcommand{\o}[0]{{\mathcal O}} 
\newcommand{\z}[0]{{\mathbb Z}}

\newcommand{\p}[0]{{\mathbb P}}

\newcommand{\qtq}[1]{\quad\mbox{#1}\quad}

\newcommand{\rank}[0]{\operatorname{rank}}

\newcommand{\rup}[1]{\lceil{#1}\rceil}

\newcommand{\onto}[0]{\twoheadrightarrow}

\newcommand{\quot}[0]{\operatorname{Quot}}

\begin{document}
\bibliographystyle{amsalpha}

\hfill\today

\title{Restrictions of stable bundles}
\author{V.~Balaji and J\'anos Koll\'ar}

\maketitle

Let $E$ be a stable vector bundle on a projective variety
$X$. The theorem of \cite{meh-ram} ensures that the restriction of
$E$ to a  general, sufficiently high degree
complete intersection curve in $X$ is again stable.
The original proof did not give  estimates for
the ``sufficiently high degree,'' but such bounds were developed
later. There are many special cases where the unstable
restrictions are fully understood, but the following three results
give the best known estimates for arbitrary bundles.
For simplicity we state them only for surfaces, which seems to be
the most subtle case.

\begin{say}[General effective results]\label{gen.results.say}
 Let $X$ be a normal, projective surface,
over an algebraicaly closed,
 $|H|$ a very ample linear system and 
 $E$ a stable vector bundle  (or reflexive sheaf) of rank $r$   on $X$. Then
\begin{enumerate}
\item \cite{fle, lan2} $E|_{D_m}$ is semistable for general  $D_m\in |mH|$
for $m\geq C_1$  where $C_1$ is roughly  $\frac14 r^2 (H^2)$.
\item  \cite{bog, lan} $E|_{D_m}$ is stable for every smooth  $D_m\in |mH|$
for  $m\geq C_2$  where $C_2$ is roughly 
 $\Delta(E):=2rc_2(E)-(r-1)c_1^2(E)$ (plus a more complicated correction term in
positive characteristic).
\item \cite{bal-par}  Assume in addition that the  characteristic is 0.
Then $E|_{D_m}$  is stable and
 has the same holonomy group \cite{bal-kol}
as $E$  for  every smooth  $D_m\in |mH|$
for $m\geq C_3$  where $C_3$ is roughly 
 $r^{r^2}\cdot\Delta(E)$.
\end{enumerate}
We refer to the original papers for more precise results and other related
bounds; see  \cite[Sec.7]{huy-leh} for an introduction.
\end{say}

For many applications, for instance for general boundedness results 
for sheaves,
any effective estimate is useful, but it would be of interest to understand
the optimal bounds on $m$ in any of the above settings.
The aim of this note is to  prove  Flenner-type theorems that
yield stability and then  suggesting a possible optimal
result along these directions.

\begin{thm}\label{main.thm.1} Let $X$ be a normal, projective surface
over an algebraicaly closed field and 
$|H|$ an ample and base point free linear system.
Let $E$ be a stable reflexive sheaf of rank $r$ on $X$ such that
$E|_C$ is semistable for general $C\in |H|$. Then
 $E|_{D_m}$ is stable for general $D_m\in |mH|$ for $m\geq \frac12{r^2}+4$.
\end{thm}

\begin{thm}\label{main.thm.2} Let $X$ be a normal, projective surface,
over an algebraicaly closed field of characteristic 0 and 
$|H|$ an ample and base point free linear system.
Let $E$ be a stable reflexive sheaf of rank $r$ on $X$ such that
$E|_C$ is stable (or polystable)  for general $C\in |H|$. Then
 $E|_{D_m}$  is stable and has the same holonomy group as $E$ 
for general $D\in |mH|$ for $m\geq 2r^2+3$.
\end{thm}

In concrete situations, we think of Theorems \ref{main.thm.1}--\ref{main.thm.2}
as building on a Flenner-type estimate. Thus first one
establishes the semistability of restrictions and then, using 
 the above results one   gets stability and the correct holonomy group
 for general $D_m\in |mH|$ with a polynomial bound on $m$.

A common 
feature of  Theorems \ref{main.thm.1}--\ref{main.thm.2}  is that they use
only {\em one}  numerical invariant, the rank of the bundle
in the estimate. We do not know how to eliminate
$(H^2)$ from the semistability bound in 
(\ref{gen.results.say}.1).

It seems easy to improve the constants in the {\em quadratic} bounds of
 Theorems \ref{main.thm.1}--\ref{main.thm.2} by more attention to
details. We did not try to optimize our  proof
since  we believe that our approach should yield a {\em linear} bound
in both cases. Even that, however, may not be optimal.
In order to call attention to how little is known,
let us pose the following (either bold or foolishly optimistic) question.

\begin{ques}\label{main.ques}  Let $X$ be a normal, projective surface,
$|H|$ a very ample linear system.
Let $E$ be a stable reflexive sheaf  on $X$.
Is  $E|_D$  stable for general $D\in |mH|$ for $m\geq 4$?
\end{ques}

\subsection*{Proof of  the Theorems}{\ }

While the claim is about the stability of $E$ when restricted to a
general, hence smooth, curve in $|mH|$, we will prove stability
for certain reducible curves.

\begin{defn} 
Let $X$ be a normal, projective surface,
$|H|$ an ample and base point free linear system (not necessarily complete).

By a {\it nodal $m$-gon with sides in $H$} we mean a  curve
$C\subset X$ that is the union of $m$ smooth members of $|H|$
and whose singularities are ordinary nodes.

We can view the space of all nodal $m$-gons 
either as a locally closed subvariety of $|mH|$ or 
 as an open subvariety of
$|H|^m$.  The latter shows that  it is an irreducible variety.
\end{defn}

\begin{prop}\label{effe.stab.prop}
Let $X$ be a normal, projective surface,
$|H|$ an ample and base point free linear system (not necessarily complete).
Fix a smooth point $x_0\in X$ and a smooth curve
$x_0\in C_0\in |H|$ such that $X$ is smooth along $C_0$.
Let $E$ be a stable reflexive sheaf of rank $r$ on $X$ such that
$E|_{C_0}$ is locally free and semistable.

 Then the  restriction of $E$ to the general nodal $m$-gon with sides in $H$
 is stable for $m\geq \frac12{r^2}+4$.
\end{prop}

Proof. Being  locally free and semistable  are open properties, hence
 $E|_{C_{\lambda}}$ is locally free and  semistable for general
$C_{\lambda}\in |H|$.

Let ${\mathcal T}=\{T_{\lambda,\mu}\}\subset |3H|$
be the space of all nodal triangles
$C_0+C_{\lambda}+C_{\mu}$ where 
$C_{\lambda}, C_{\mu}\in |H|$,  $X$ is smooth along $C_{\lambda}, C_{\mu}$
 and the restrictions 
$E|_{C_{\lambda}},E|_{C_{\mu}}$
 are   locally free and semistable.
Let $C_{\mathcal T}\subset X\times {\mathcal T}$
be the universal curve and
$E_{\mathcal T}\to C_{\mathcal T}$ the pull-back of $E$ to 
$C_{\mathcal T}$. We can also think of $E_{\mathcal T}$ as the 
universal vector bundle
whose restriction to $T_{\lambda,\mu}$ is
 $E_{\lambda,\mu}:=E|_{T_{\lambda,\mu}}$.

Let  $D_{\mathcal T}\subset \quot\bigl(E_{\mathcal T}\bigr)$
denote the subscheme parametrizing torsion free quotients 
$q_{\lambda,\mu}:E_{\lambda,\mu}\onto F$
with the same 
slope as $E_{\lambda,\mu}$. 
Its fiber over $T_{\lambda,\mu}$ is denoted by 
 $D_{\lambda,\mu}\subset \quot\bigl(E_{\lambda,\mu}\bigr)$.

Note that $E_{\lambda,\mu}$ is  semistable by
(\ref{reducible.semistab}) and 
each $q_{\lambda,\mu}:E_{\lambda,\mu}\onto F$ induces 3 quotients
$$
E|_{C_0}\to F_{0}, \quad 
E|_{C_{\lambda}}\to F_{\lambda}\qtq{and} E|_{C_{\mu}}\to F_{\mu}.
$$
Each of these has the same slope as $E|_{C_0}$
and at each node  $p\in C_0+C_{\lambda}+C_{\mu}$
the two branches give the same quotient
$E_p\to F_p$.
In particular, 
any such $q_{\lambda,\mu}:E_{\lambda,\mu}\onto F$ is uniquely determined
by 
$$
q_{\lambda,\mu}\otimes k(x_0): E_{x_0}\to F_{x_0}.
$$
Therefore  we can think of
$D_{\lambda,\mu}$ as a subscheme of
 $\quot\bigl(E_{x_0}\bigr)$,
which  is a union of Grassmannians  of quotients of the vector space $E_{x_0}$.
Thus 
$$
\dim D_{\lambda,\mu}\leq \dim \quot\bigl(E_{x_0}\bigr)\leq \tfrac14 r^2.
$$
Note that since each $E_{\lambda,\mu}$ is  semistable,
$D_{\mathcal T}\to {\mathcal T}$ is proper.

We apply (\ref{irred.lem}) to 
$D_{\mathcal T}\subset {\mathcal T}\times \quot\bigl(E_{x_0}\bigr)$
and first we show that the alternative (\ref{irred.lem}.1) is impossible if  
$E$ is  stable.

If  (\ref{irred.lem}.1) holds then $D_{\mathcal T}\to {\mathcal T}$
has a constant section. Equivalently,
there is a quotient
$$
q_{\mathcal T}: E_{\mathcal T}\to F_{\mathcal T},
$$
such that, for every $T_{\lambda,\mu}\in {\mathcal T}$,
its restriction gives a   quotient  (which has the same 
slope as $E_{\lambda,\mu}$) 
$$
q_{\lambda,\mu}:E_{\lambda,\mu}\onto F_{\lambda,\mu}
\qtq{such that} q_{\lambda,\mu}\otimes k(x_0)=q_0.
$$
We use these to construct a quotient sheaf  $E\to F$
whose pull-back to $C_{\mathcal T}$ is $F_{\mathcal T}$.
This will then contradict the
stability of $E$.

To construct $E\to F$, 
pick a general point $x\in X$ and let
$C_{\lambda},C_{\mu}\in |H|$ be  general smooth curves through $x$.
Set $q_x:= q_{\lambda,\mu}\otimes k(x)$.

Note that $E|_{C_0+C_{\lambda}}$ is semistable, hence
$q_{\lambda,\mu}|_{C_0+C_{\lambda}}$ is uniquely determined by
 $E|_{C_0+C_{\lambda}}$ and by $q_0$. Thus
$q_x$ does not depend on the choice of $C_{\mu}$.
Similarly, it also  does not depend on the choice of $C_{\lambda}$.
Thus, as the notation suggests, 
$q_x$ is independent of the choices of $C_{\lambda},C_{\mu}$.
Thus we get a well defined quotient $q:E\to F$ 
such that $F|_{C_{\lambda}}=F_{\lambda,\mu}|_{C_{\lambda}}$
for general $C_{\lambda}\in |H|$.
This shows that
$E$ is not stable, a contradiction.
Therefore the alternative (\ref{irred.lem}.2) must hold. 

Fix  $m_0\geq \frac{r^2}{4}+1$ and pick general
pairs  $C_{\lambda_i},C_{\mu_i}\in |H|$.
We claim that  
the restriction of $E$ to the $(2m_0+1)$-gon
$C_{\Sigma}:=C_0+\sum_{i=1}^m \bigl(C_{\lambda_i}+C_{\mu_i}\bigr)$ is stable.

Assume to the contrary that  there is a
quotient 
$q_{\Sigma}:E|_{C_{\Sigma}}\to F_{\Sigma}$ which has the same slope as
$E|_{C_{\Sigma}}$. 
Let $q_0:E_{x_0}\to F_{x_0}$ be the induced quotient on the fiber over $x_0$.
We can restrict $q_{\Sigma}$ to
$$
q_i: E|_{C_0+C_{\lambda_i}+C_{\mu_i}}\to F_i\qtq{for $ i=1,\dots, m$.}
$$
Each $q_i$ gives a point in
$D_{\lambda_i,\mu_i}$ and
$[q_0]\in \cap_{i=1}^m D_{\lambda_i,\mu_i}$.
This contradicts (\ref{irred.lem}.2), hence $E|_{C_{\Sigma}}$ is stable.

The smallest value of $m_0$ we can take is
 $ \rup{\frac14{r^2}}+1$, which gives that
$E$ restricted to the general $m$-gon is stable
for $m\geq 2\rup{\frac14{r^2}}+3$. This holds if 
 $m\geq \frac12{r^2}+4$.
 \qed

\begin{say}[Proof of  (\ref{main.thm.1})]
By  (\ref{effe.stab.prop}),  $E$ restricted to the  $m$-gon is stable.
Since stability is an open condition, 
 $E$ restricted to the general member of
$|mH|$ is also stable. \qed
\end{say}

\begin{lem} \label{irred.lem}
Let $U$ be an irreducible variety, $V$ any variety
and $Z\subset U\times V$ a closed subscheme. Then
\begin{enumerate}
\item either $U\times \{v\}\subset Z$ for some $v\in V$, 
\item or  $\cap_{i=1}^m Z_{u_i}=\emptyset$
for $m>\dim V$ and 
 general $u_1,\dots, u_m\in U$.
\end{enumerate}
\end{lem}

Proof. Assume that (1) fails. By induction
we show that $\dim \cap_{i=1}^r Z_{u_i}\leq \dim V-r$ for
 $r\leq \dim V+1$  and general $u_1,\dots, u_r\in U$.

This is clear for $r=0$. To go from $r$ to $r+1$,
note that  none of the irreducible components
of $\cap_{i=1}^r Z_{u_i}$ is contained in every $Z_u$.
Thus, a general $Z_{u_{r+1}}$ contains none of them, hence
$$
\dim \cap_{i=1}^{r+1} Z_{u_i}< \dim \cap_{i=1}^r Z_{u_i}. \qed
$$

\begin{prop}\label{reducible.semistab}
 Let $E$ be a stable (resp.\ semistable)
bundle on $X$ and let $C_1,\dots,C_r \in |H|$ be smooth curves such that the
$E|_{C_i}$ are stable (resp.\ semistable). Then
 $E|_{C_1+\cdots+C_r}$ is also stable (resp.\ semistable). 
\end{prop}

Proof.  This follows essentially from \cite[Prop.1.2]{bigas}. The only thing
 to observe  is that, since the curves all lie on the surface $X$, the
 weights that \cite[Sec.7]{sesh} associates to torsion-free sheaves on 
$C_1+\cdots+C_r$  for the purposes of defining semi-stability are all equal.
 From this, one
 observes that the inequality  \cite[1.1]{bigas} is satisfied in 
 our situation.  \qed

\begin{rem} It is easy to modify the above results
to get stability not only for  general $C\in |mH|$ but also
for general   $C\in |mH|$ passing through some preassigned points.

Let $F$ be stable on $X$ 
and fix $x_1,\dots, x_r\in X$ such that $F$ is locally free at these points.
 Let $p:Y \to X$ be the blow-up of $X$ along the points  $x_1,\dots, x_r$
with  exceptional divisor $D$.

Set $E = p^{*}F$. Then, as in  \cite[Prop.3.4]{buch},
for some $n \gg 1$, $E$ is stable on $Y$ 
with respect to the polarization
 $H_n:= np^{*}(H) - D$. Fix this $n$ and set $H_Y:=H_n$. 

Now apply the above results on stable bundles to get
 $m$ such for  $E|_{C_0}$ is stable for
sufficiently general $C_0 \in |mH_Y|$. 

Consider $E|_{C_0 + D}$. Since $E|_D$ is trivial, it is semistable with 
respect to  $H_Y$ and 
since $E|_{C_0}$ is stable, by \cite[Prop.1.2]{bigas} it follows that 
$E|_{C_0 + D}$ is actually stable.
Being stable is an open condition, hence 
 we see that  $E|_{C_1}$ is  stable for a general member 
$C_1 \in |mH_Y + D| = |mp^*H - (m-1)D|$. Proceeding 
the same way, we get eventually that $E|_{C_{{m-1}}}$ is stable for a general
 member
 $C_{m-1} \in |mH_Y + (m-1)D| =|mp^*H - D|$. 
Set $C:=p(C_{m-1})$ and note that in fact $C\cong C_{m-1}$. Thus 
 $V|_{C}$ is stable and 
$C$ is a smooth member of $|mH|$ passing through the points $x_1,\dots, x_r$.  
  \end{rem}

\begin{say}[Proof of  Theorem \ref{main.thm.2}]\label{pf.of.2}

The arguments are quite similar to the ones
used to show Theorem \ref{main.thm.1}, hence we only
outline them.

As usual (see, for instance,  \cite[\S 3]{bal-kol}), 
by passing to a finite cover of $X$ if
necessary, we may assume that $\det E$ is trivial.

Pick a general point $x_0\in X$ and a curve
$x_0\in C_0\in |H|$ such that 
$E_{C_0}$ is  locally free and polystable.
By assumption 
 $E|_{C_{\lambda}}$ is locally free and  polystable for general
$C_{\lambda}\in |H|$.
Since $\deg E|_{C_{\lambda}}=0$, by \cite{nar-ses}, 
one can also obtain 
 $E|_{C_{\lambda}}$ from a unitary representaion of
the fundamental group of $C_{\lambda}$. In particular,
there is a well defined notion of parallel transport
along any path in $C_{\lambda}$ or in any $m$-gon
$\cup_iC_{\lambda_i}$ if $E|_{C_{\lambda_i}}$ is locally free and  
polystable for every $i$.

Let ${\rm Hol}_{x_0}(E)\subset GL(E_{x_0})$
denote the holonomy group \cite{bal-kol} and
 ${\rm Hol}^{\circ}_{x_0}(E)\subset {\rm Hol}_{x_0}(E)$
 its identity component.
By \cite[40]{bal-kol}, 
$\pi_1(X, x_0)\to {\rm Hol}_{x_0}(E)/{\rm Hol}^{\circ}_{x_0}(E)$
is surjective. 
By the Lefschetz theorem, $\pi_1(C_0, x_0)\to \pi_1(X, x_0)$ is surjective,
which implies that
$$
{\rm Hol}_{x_0}(E|_{C_0})\to {\rm Hol}_{x_0}(E)/{\rm Hol}^{\circ}_{x_0}(E)
\qtq{is surjective.}
\eqno{(\ref{pf.of.2}.1)}
$$
Therefore, using
 \cite[40]{bal-kol}, by passing to a suitable finite \'etale cover of $X$
we may assume that
 ${\rm Hol}_{x_0}(E)$ is connected.

Let ${\mathcal T}=\{T_{\lambda,\mu}\}\subset |3H|$
be the space of all nodal triangles
$C_0+C_{\lambda}+C_{\mu}$ where 
$C_{\lambda}, C_{\mu}\in |H|$,  $X$ is smooth along $C_{\lambda}, C_{\mu}$
 and the restrictions 
$E|_{C_{\lambda}},E|_{C_{\mu}}$
 are   locally free and polystable.

Using the parallel transport along the 3 components, 
we get a (possibly non-unitary) representation
$$
\rho_{\lambda,\mu}: \pi_1\bigl(C_{\lambda,\mu}, x_0\bigr)\to 
GL\bigl(E_{x_0}\bigr).
$$

The images of all these  representations
generate a subgroup $H_{\mathcal T}$ of the holonomy group 
${\rm Hol}_{x_0}(E)$. 

First we claim  that $H_{\mathcal T}={\rm Hol}_{x_0}(E)$.
Assume to the contrary that 
$H_{\mathcal T}\subsetneq {\rm Hol}_{x_0}(E)$.
Then there is a tensor power $E_{x_0}^{\otimes m}$
and a vector $w\in E_{x_0}^{\otimes m}$ that is
$H_{\mathcal T}$-invariant but not ${\rm Hol}_{x_0}(E)$-invariant.
Thus, for each triangle $C_{\lambda,\mu}$
we get a flat section
$$
w_{\lambda,\mu}\in H^0\Bigl(C_{\lambda,\mu}, \bigl(E_{\lambda,\mu}\bigr)^{\otimes m}\Bigr).
$$
As in the proof of Theorem \ref{main.thm.1}
we see that for every
$x\in C_{\lambda}\cap C_{\mu}$, the fiber
$w_{\lambda,\mu}(x)$ depends only on $x$ but not on
$C_{\lambda}$ and $C_{\mu}$. Thus we get a well defined global section
$$
w_X\in H^0\bigl(X\setminus (\mbox{finite set}),  E^{\otimes m}\bigr)
$$
which then extends to a global section of (the reflexive hull of)  $E^{\otimes m}$.
Thus $w=w_X(x)$ is  ${\rm Hol}_{x_0}(E)$-invariant, a contradiction.
This proves that $H_{\mathcal T}={\rm Hol}_{x_0}(E)$.

Continuing with the method of  Theorem \ref{main.thm.1}
would give a bound that depends on $r$ and on $m$ above.
In many important cases, for instance when
${\rm Hol}_{x_0}(E)=SL(E_{x_0})$, one can choose $m=2$ \cite[Prop.5]{bal-kol}.
 However, even this would give a degree 4 bound in $r$.
In general, it is not known how to bound $m$ effectively.

Thus, instead of trying to control the quot-scheme as in
the proof  of  Theorem \ref{main.thm.1},
we control the size of the holonomy group on $m$-gons
 using (\ref{generate.lem}).

Choose $r^2+1$ general pairs
$(\lambda_i, \mu_i)$. Then, by   (\ref{generate.lem}),
the images of $\rho_{\lambda_i,\mu_i}$ for $i=1,\dots, r^2+1$
 generate ${\rm Hol}_{x_0}(E)$.

This implies that the restriction of $E$ to the general
$(2r^2+3)$-gon is stable and has holonomy group
${\rm Hol}_{x_0}(E)$. 
By the lower semicontinuity of the  holonomy groups
 \cite[\S 1]{bal-kol},
 the same holds for a general smooth curve in
$\bigl|(2r^2+3)H\bigr|$.\qed

\begin{lem} \label{generate.lem}
Let $G\subset GL(n,\c)$
 be a connected algebraic group of dimension $d$ over $\c$
and $S\subset G$ a connected (in the Euclidean topology) subset
that generates a Zariski dense subgroup of $G$.
Then there are $d+1$ elements $s_0,\dots, s_d\in S$ that
 generate a Zariski dense subgroup of $G$.
\end{lem}

Proof. For $0\leq r\leq d$ we use  induction to find $s_0,\dots, s_r\in S$
such that the  Zariski  closure of
$\langle s_0,\dots, s_r\rangle$ has dimension at least $r$.
 This is clear for $r=0$.

To start with, fix any $s_0\in S$ and  consider
$ss_0^{-1}$ as a function $S\to G$. It sends $s_0$ to the identity.
If the eigenvalues of $ss_0^{-1}$ 
are not constant near $s_0$,
then for very general $s_1\in S$, at least one of the 
eigenvalues of $s_1s_0^{-1}$ is not a root of unity.
Then  $s_1s_0^{-1}$ has infinite order, hence the
Zariski  closure of $\langle s_1s_0^{-1}\rangle$ has positive dimension.
If the eigenvalues of $ss_0^{-1}$ are  constant near $s_0$,
then all the  eigenvalues of $s_1s_0^{-1}$  equal 1.
Thus $s_1s_0^{-1}$ has infinite order, unless $s_1=s_0$.

Now to the inductive step.
Let $H_r\subset G$ denote the  Zariski  closure of
$\langle s_0,\dots, s_r\rangle$ and $H^{\circ}_r\subset H_r$ its
 identity component.
By assumption, $\dim H^{\circ}_r\geq r$.

If $H^{\circ}_r$ is a normal subgroup of $G$, then the above argument applies to
$G/ H^{\circ}_r$. Pick $s_{0r}\in S\cap H^{\circ}_r$ 
such that no open neighborhood of $s_{0r}\in U\subset S$
is contained in $H^{\circ}_r$.
We obtain that  the  Zariski  closure of
$\langle s_{r+1}s_{0r}^{-1}\rangle$ is a positive dimensional
subgroup of $G/ H^{\circ}_r$. Thus
the  Zariski  closure of
$\langle s_0,\dots, s_r, s_{r+1}s_{0r}^{-1}\rangle$ has dimension at 
least $r+1$.
Since $s_{0r}$ is in the Zariski  closure of
$\langle s_0,\dots, s_r\rangle$,
we can replace $s_{r+1}s_{0r}^{-1}$ by $s_{r+1}$
without changing the Zariski closure.

 If $H^{\circ}_r$ is not a normal subgroup of $G$, then pick an
$s_{r+1}$ that is not contained in the normalizer of  $H^{\circ}_r$.
Then $H^{\circ}_r$ is not a normal subgroup of
 the  Zariski  closure of $\langle s_0,\dots, s_{r+1}\rangle$.
The  identity component is always a  normal subgroup, thus
the  identity component  of
 the  Zariski  closure of $\langle s_0,\dots, s_{r+1} \rangle$
is strictly larger than $H^{\circ}_r$.\qed
\end{say}

\begin{rem} (1) The connectedness of $S$ is essential in
(\ref{generate.lem}). For instance, all the
roots of unity 
 generate a Zariski dense subgroup of $\c^*$, but any finite
subset of them  generates a finite subgroup.

(2) It is easy to see that 2 very general elements of a
connected, reductive, algebraic group  generate a Zariski dense subgroup.
Indeed,  the Zariski closure of the subgroup generated by 
a very general semisimple element $g_1$
is a maximal torus. The maximal torus acts on the Lie algebra of $G$
with 1-dimensional eigenspaces (except on the  Lie algebra of the torus),
 hence only finitely many
connected subgroups contain any given maximal torus.
Pick any $g_2\in G$ not in  the normalizer of any of these
subgroups that are not normal in $G$. Then $\langle g_1, g_2 \rangle$
is a  Zariski dense subgroup of $G$.

(3) Probably a small case analysis would improve the bound
$\dim G+1$ in (\ref{generate.lem}) to $\dim G$ which is the optimal result
for $G=\c^d$, where $d-1$ elements always generate a smaller
dimensional subgroup.
A very general element of  $\bigl(\c^*\bigr)^d$
generates a  Zariski dense subgroup, but if we take
$S\subset \bigl(\c^*\bigr)^d$ to be the
union of the ``coordinate axes''
$(1,\dots, 1, *, 1,\dots, 1)$ then again no
$(d-1)$-element subset of $S$ generates   $\bigl(\c^*\bigr)^d$.

We believe, however, that one can do much better for
reductive groups, especially if $S\subset G$ is an
irreducible real algebraic subset. Here the worst example we
 know is the following.

(4) The set of all reflections generate the orthogonal group
$O(d)$ but  $d-1$ reflections always have a common fixed
vector, hence they generate a smaller dimensional subgroup.
(The orthogonal group is not connected, so it may be better to
work with the orthogonal similitudes and with
scalars times reflections.)
\end{rem}

\begin{ques} \label{generate.ques}
Let $G$ be a connected, reductive algebraic group of rank  $r$ over $\c$
and $S\subset G$ an irreducible, real, semialgebraic subset
that generates a Zariski dense subgroup of $G$.
Is it true that  $2r$ very general elements $s_1,\dots, s_{2r}\in S$ 
 generate a Zariski dense subgroup of $G$.
\end{ques}

\subsection*{Remarks on Question \ref{main.ques}}{\ }

More generally, one can investigate the following.

\begin{ques}\label{main.gen.ques}  Let $X$ be a smooth, projective surface and
$|H|$ an ample and base point free linear system. 
Under what conditions on $(X,|H|)$ can one guarantee that for every 
 stable vector bundle $E$ on $X$, the restriction
  $E|_C$  is stable for general $C\in |H|$?
\end{ques}

We know very few examples of pairs  $(X,|H|)$ where
stability of restrictions fails.
One such case is when  a general $C\in |H|$ is rational or elliptic.
This holds, among others,  for $\bigl(\p^2, |\o_{\p^2}(1)|\bigr)$,
 $\bigl(\p^2, |\o_{\p^2}(2)|\bigr)$ and $\bigl(\p^2, |\o_{\p^2}(3)|\bigr)$.

On a rational curve every stable bundle has rank 1 and on an
elliptic curve every stable bundle with $c_1(E)=0$ has rank 1.
Thus if $\rank E\geq 2$ and  $c_1(E)=0$ then $E|_C$ is never stable.

One can get more complicated examples out of these.
Take any surface $X$ and a general finite morphism
$\pi:X\to \p^2$. Set $|H|:=\pi^*|\o_{\p^2}(3)|$.
Note that $|H|$ is ample but usually neither very ample nor complete.
There are, however, many examples, for instance double covers
whose branch locus has degree $\geq 8$, where
the pulled-back $|H|$ is a complete linear system whose general member
is a smooth curve of high genus. Nonetheless, if
$E$ is the pull-back of a vector bundle from $\p^2$, then
the  restrictions
$E|_C$ are not stable for $C\in |H|$.

These examples all satisfy $\dim |H|\leq 9$, but
pulling back  $|\o_{\p^1\times \p^1}(1,a)|$ gives similar examples
where both $\dim |H|$ and the genus of the general $C\in |H|$
are arbitrarily high.

These types are the obvious  examples where
general restrictions are not stable. We do not know any other.

Let us next turn to a heuristic argument that suggested
 Question \ref{main.ques} to us.
We focus on the holonomy groups and propose the following
variant.

\begin{ques}\label{main.hol.ques}  Let $X$ be a smooth projective surface and
$|H|$ a very ample linear system.
Let $E$ be a stable vector bundle  on $X$.
Is  $E|_D$  stable for general $D\in |mH|$ for $m\geq 4$
and with the same holonomy group as $E$?
\end{ques}

As we saw in (\ref{pf.of.2}.1),  the discrete part of the holonomy
${\rm Hol}(E)/{\rm Hol}^{\circ}(E)$ never causes problems
in (\ref{main.hol.ques}). Thus, by
 \cite[40]{bal-kol}, we can  focus on the case
when ${\rm Hol}(E)$ is connected and $\det E\cong \o_X$.

For $m\gg 1$  take a general  $x\in C\in |mH|$ such that
 $E|_C$ is stable and  with the same holonomy group as $E$.
Thus we get a holonomy representation
$$
\rho_C:\pi_1(C)\to {\rm Hol}(E)
$$
whose image is Zariski dense.

Although not supported by any evidence, one can hope
that in our situation $C$ can be written as a connected sum 
$C=C_2\# C'$ where the
genus of $C_2$ is 2 and
$$
\rho_2:\pi_1(C_2)\to {\rm Hol}(E)
\qtq{still has Zariski dense image.}
$$
It is then another entirely uncorroborated belief
that this $C_2$ can be realized by vanishing cycles
as the curve acquires an ordinary 4-fold point.
Eventually, this may lead to an approximation of
$\rho_2:\pi_1(C_2)\to {\rm Hol}(E)$
by some $\rho_t:\pi_1(C_t)\to {\rm Hol}(E)$
where $C_t\in |4H|$ is a family of curves
whose limit also has  an ordinary 4-fold point.

We stress that for the moment all this is just wishful thinking.
We, however, feel that this approach  raises many interesting questions that --
even if  Questions \ref{main.ques} and \ref{main.hol.ques}
turn out to be utterly misguided --
could lead to a much improved understanding of
stable bundles and their restrictions.

\begin{ack} We thank  I.~Coskun, A.~Langer, 
 A.J.~Parameswaran and C.S.~Seshadri 
 for useful  comments and questions.
Partially support for VB was provided by the J.C.~Bose research
grant.
Partial financial support for JK  was provided by  the NSF under grant number 
DMS-0758275.
\end{ack}

\bibliography{refs}

\providecommand{\bysame}{\leavevmode\hbox to3em{\hrulefill}\thinspace}
\providecommand{\MR}{\relax\ifhmode\unskip\space\fi MR }
\providecommand{\MRhref}[2]{%
  \href{http://www.ams.org/mathscinet-getitem?mr=#1}{#2}
}
\providecommand{\href}[2]{#2}
\begin{thebibliography}{TiB95}

\bibitem[BK08]{bal-kol}
V.~Balaji and J{\'a}nos Koll{\'a}r, \emph{Holonomy groups of stable vector
  bundles}, Publ. Res. Inst. Math. Sci. \textbf{44} (2008), no.~2, 183--211.
  \MR{2426347 (2010c:14044)}

\bibitem[Bog94]{bog}
F.~A. Bogomolov, \emph{Stable vector bundles on projective surfaces}, Mat. Sb.
  \textbf{185} (1994), no.~4, 3--26. \MR{1272185 (95j:14056)}

\bibitem[BP11]{bal-par}
V.~{Balaji} and A.~J. {Parameswaran}, \emph{An analogue of the
  {N}arasimhan-{S}eshadri theorem in higher dimensions and some applications},
  Journal of Topology (to appear, math.AG:0809.376) (2011).

\bibitem[Buc00]{buch}
Nicholas~P. Buchdahl, \emph{Blowups and gauge fields}, Pacific J. Math.
  \textbf{196} (2000), no.~1, 69--111. \MR{1797236 (2001m:32038)}

\bibitem[Fle84]{fle}
Hubert Flenner, \emph{Restrictions of semistable bundles on projective
  varieties}, Comment. Math. Helv. \textbf{59} (1984), no.~4, 635--650.
  \MR{780080 (86m:14014)}

\bibitem[HL97]{huy-leh}
Daniel Huybrechts and Manfred Lehn, \emph{The geometry of moduli spaces of
  sheaves}, Aspects of Mathematics, E31, Friedr. Vieweg \& Sohn, Braunschweig,
  1997. \MR{1450870 (98g:14012)}

\bibitem[Lan04]{lan}
Adrian Langer, \emph{Semistable sheaves in positive characteristic}, Ann. of
  Math. (2) \textbf{159} (2004), no.~1, 251--276. \MR{2051393 (2005c:14021)}

\bibitem[Lan10]{lan2}
\bysame, \emph{A note on restriction theorems for semistable sheaves}, Math.
  Res. Lett. \textbf{17} (2010), no.~5, 823--832. \MR{2727611}

\bibitem[MR84]{meh-ram}
V.~B. Mehta and A.~Ramanathan, \emph{Restriction of stable sheaves and
  representations of the fundamental group}, Invent. Math. \textbf{77} (1984),
  163--172.

\bibitem[NS65]{nar-ses}
M.~S. Narasimhan and C.~S. Seshadri, \emph{Stable and unitary vector bundles on
  a compact {R}iemann surface}, Ann. of Math. (2) \textbf{82} (1965), 540--567.
  \MR{MR0184252 (32 \#1725)}

\bibitem[Ses82]{sesh}
C.~S. Seshadri, \emph{Fibr\'es vectoriels sur les courbes alg\'ebriques},
  Ast\'erisque, vol.~96, Soci\'et\'e Math\'ematique de France, Paris, 1982,
  Notes written by J.-M. Drezet from a course at the {\'E}cole Normale
  Sup{\'e}rieure, June 1980. \MR{699278 (85b:14023)}

\bibitem[TiB95]{bigas}
Montserrat Teixidor~i Bigas, \emph{Moduli spaces of vector bundles on reducible
  curves}, Amer. J. Math. \textbf{117} (1995), no.~1, 125--139. \MR{1314460
  (96e:14014)}

\end{thebibliography}

\bigskip

\noindent Chennai Math.\ Inst.\ SIPCOT IT Park, Siruseri-603103, India

\begin{verbatim}balaji@cmi.ac.in\end{verbatim}

\bigskip

\noindent Princeton University, Princeton NJ 08544-1000, USA

\begin{verbatim}kollar@math.princeton.edu\end{verbatim}

\end{document}